\documentclass{article}

\usepackage{algorithm}
\usepackage{algorithmic}
\usepackage{authblk}
\usepackage{biblatex}
\usepackage[hscale=0.7,vscale=0.8]{geometry}
\usepackage{graphicx}
\usepackage{adjustbox}
\usepackage{amsmath}
\usepackage{hyperref}
\usepackage{threeparttable}
\usepackage{tikz}
\usepackage[section]{placeins}
\usetikzlibrary{arrows.meta}
\usetikzlibrary{fit}
\allowdisplaybreaks
\definecolor{darkgreen}{RGB}{0,150,0}
\title{Robust Railway Network Design based on Strategic Timetables}

\bibliography{literature}	

\begin{document}

\author[1]{Tim Sander \thanks{\href{mailto:tim.sander@tu-dresden.de}{tim.sander@tu-dresden.de}} }
\author[2]{Nadine Friesen}
\author[1]{Karl Nachtigall}
\author[2]{Nils Nießen}

\affil[1]{Chair of Traffic Flow Science, TU Dresden, Dresden, Germany}
\affil[2]{Institute of Transport Science, RWTH Aachen, Aachen, Germany}

\maketitle

\begin{abstract}
Using strategic timetables as input for railway network design has become increasingly popular among western European railway infrastructure operators. Although both railway timetabling and railway network design on their own are well covered by academic research, there is still a gap in the literature concerning timetable-based network design. Therefore, we propose a mixed-integer linear program to design railway infrastructure so that the demand derived from a strategic timetable can be satisfied with minimal infrastructure costs. The demand is given by a list of trains, each featuring start and destination nodes as well as time bounds and a set of frequency and transfer constraints that capture the strategic timetable's main characteristics. During the optimization, the solver decides which railway lines need to be built or expanded and whether travel or headway times must be shortened to meet the demand. Since strategic timetables are subject to uncertainty, we expand the optimization model to a robust version. Uncertain timetables are modelled as discrete scenarios, while uncertain freight train demand is modelled using optional trains, which can be inserted into the resulting timetable if they do not require additional infrastructure. We present computational results for both the deterministic and the robust case and give an outlook on further research.

\end{abstract}
	
\section{Introduction}\label{sec:Introduction}
Many western European railway infrastructure operators, including German Deutsche Bahn (DB) Netz AG and Swiss Schweizerische Bundesbahn (SBB), plan expansions of their railway networks according to a long-term, strategic timetable. These timetables are typically nationwide and follow the principles of the integrated timetable: trains travel in fixed frequencies, and connections in important stations are optimized. However, creating these strategic timetables relies heavily on the planner's expertise and involves a lot of manual, albeit software-supported work. The same is true for the following identification of infrastructure expansions which are necessary to operate the strategic timetable. To support this process through mathematical optimization, we propose an optimization model for network design based on strategic timetables. It determines the number of tracks for each section and necessary improvements of travel and headway times to satisfy the demand derived from a given strategic timetable, which includes individual trains, line frequencies, and ensured connections. Since strategic timetables are constructed many years in advance, they are subject to uncertainty induced by varying demand prognoses or political circumstances. Our network design model accounts for these uncertainties by allowing the calculation of networks suitable for a timetable family consisting of several timetable scenarios and varying demand within one scenario. 
	
The remainder of this paper is structured as follows: chapter~\ref{sec:Literature} summarizes the current state of the literature. The following chapter~\ref{sec:Methodology} describes how we model infrastructure (section~\ref{sec:ModellingInfrastructure}), input timetables (section~\ref{sec:ModellingTimetables}) and uncertainties (section~\ref{sec:ModellingUncertainty}). The deterministic optimization model is introduced in section~\ref{sec:DeterministicModel}, followed by the robust expansion in section~\ref{sec:RobustModel}. Chapter~\ref{sec:CaseStudy} describes the case study and provides computational results. Finally, we provide a summary and an outlook on further research in chapter~\ref{sec:SummaryOutlook}.
	
\section{Literature Review}\label{sec:Literature}
	
Railway network design based on strategic timetables touches on several research topics. This chapter summarizes recent contributions to the topics of timetable-based railway infrastructure development, railway network design, strategic and tactical timetabling and service network design with timing consideration.
	
The first time that the concept of using a strategic timetable as the base for designing railway infrastructure was utilized was in Switzerland in the 1980s when the first ideas for the concept \textit{Bahn 2000} were conceived. A nationwide timetable was designed and used as a reference to plan infrastructure expansions and new high-speed lines. This strategic timetable featured a line concept with trains running in fixed intervals of 30 or 60 minutes and integrated connections in important stations. See~\cite{Kraeuchi_Bahn_2000} for a detailed, although not scientific, description.
This concept has been adapted in Germany and has been described in several articles, including~\cite{Weigand.2012} or~\cite{Korner.2012}. In these articles, employees of the largest German railway infrastructure manager, DB Netz AG, describe the concept of the so-called long-term timetable and its role in railway capacity management. Also, the idea of using the long-term timetable as input for strategic network design is further illustrated. However, these articles only describe general concepts and do not provide optimization models for either timetabling or network design. In~\cite{Heppe2019}, Heppe and Weigand describe the current process of long-term timetable planning and infrastructure planning but without elaborating on methods or providing an optimization model.
With the so-called \textit{Deutschlandtakt}, a nation-wide integrated timetable for Germany has been developed~\cite{SMAundPartnerAG.Deutschlandtakt.2021} and resulting infrastructure demands have been identified~\cite{SMAundPartnerAG.Infrastruktur.2021}. However, the strategic timetable and the list of infrastructure were created using traditional planning methods and without optimization models.
	
While the general concept of timetable-based network design for railway infrastructure is well described, there appears to be a lack of optimization models for this use case. However, several papers describe optimization models for timetable-independent railway network design. Spönemann describes in his PhD thesis~\cite{Sponemann.2013} an approach which models the railway infrastructure on a macroscopic level as a multigraph, where each of the parallel arcs represents one level of infrastructure expansion. He uses minimal headway times to estimate the capacity consumption on railway lines. Based on a timetable-independent demand estimation, the number of tracks per direction and the distance between overtaking opportunities necessary to satisfy the demand on each track section are determined. A different approach is shown in \cite{Seyedvakili.2018}, where the authors developed an optimization model which chooses the most efficient infrastructure expansions out of a list of proposed measures. They further expanded their model in~\cite{Seyedvakili.2020} by considering several planning periods.
	
Railway network design has also been addressed by~\cite{Laporte.2010} and ~\cite{Laporte.2011} while focussing on heavy-rail rapid transit networks. They also take link failures into account and provide a game-theoretic solution approach. A similar model is provided by~\cite{GarciaArchilla.2013}, where the usage of a GRASP algorithm is proposed to solve the problem. These approaches focus on public transit within a city or metropolitan area and are closely related to the public transit network design problem. This class of models decides where to build stations and how to connect them while often considering line planning, passenger routing and robustness towards disruptions. A recent survey on this topic can be found in~\cite{Laporte.2019}. The public transit network design problems usually use a few simplifications, e.g. assuming homogeneous rolling stock and a completely double-tracked network, and are thus not suitable for our research focus.
	
However, despite this large number of publications, there exist, to the best of the author's knowledge, only two papers, written by the same authors, that explicitly consider the input of a strategic timetable for network design:~\cite{Schobel.2012} and~\cite{Schobel.2013}. Both publications focus on microscopic infrastructure modelling. In~\cite{Schobel.2012}, a graphical solution approach to determine the location of double-track sections on a single-track line is introduced and refined using a microscopic simulation tool. The authors describe an optimization model to solve the network design problem in~\cite{Schobel.2013}. However, no working solution method is provided. Also, the microscopic infrastructure model limits the size of instances that can be optimized in reasonable time severely.
	
Besides the application to railways and public transit, network design problems, in general, are a vast field of research which has seen many publications in the last decades. They were first introduced and detailed in~\cite{Wong.1976} and~\cite{Magnanti.1984}. A more recent collection of problem formulations, solution algorithms and applications to transportation and logistics can be found in~\cite{Crainic.2021}.
	
Another large field of research concerning our own studies is railway timetabling. Many different model formulations and solution approaches are developed. Several models also integrate timetabling with one of the other planning steps, like crew and vehicle scheduling or line planning. However, most timetabling research focuses on the tactical planning stage and does not use a strategic perspective. For recent reviews of railway timetabling, please refer to chapters five and six of~\cite{Borndoerfer.2018} or to~\cite{Caimi.2017}. Only recently, models focusing explicitly on strategic timetables have been introduced by~\cite{Polinder.2021},~\cite{Graaf.2021} and~\cite{coviello.2021}. These approaches focus on the passenger perspective and either assume a fixed infrastructure, including known projects that are not yet finished but will be available once the strategic timetable becomes operational, or relax infrastructure constraints. The model of~\cite{Polinder.2021}, which relaxes infrastructure constraints, has been further developed in~\cite{Polinder.2022} by using a fixed infrastructure model on which the passenger-oriented timetabling problem is described and solved. In~\cite{Sartor.2023}, a quasi-periodic timetabling model is used to evaluate different infrastructure expansion scenarios.
	
When researching the integration of timetable-based constraints into network design and thus including a temporal perspective into the model, one comes across service network design (SND) and service network scheduling problems. These problems often consider consolidation-based transportation companies, characterized by grouping individual passengers of freight loads for different customers into one trip and one vehicle~\cite{Crainic.2021}. They optimize the operations of a transport company while fulfilling customer demand and minimizing spent resources. An overview of service network design problems and some of their applications can be found in chapters 12 to 17 of~\cite{Crainic.2021}. SND problems which consider the timing and scheduling of trips often use time-expanded networks. Here, the network of nodes and arcs is duplicated for each time step to model timed events. However, choosing the right time step comes at a trade-off between solution quality and computation times. Short discretization intervals can lead to computationally intractable models while larger intervals reduce the solution quality and accuracy. This price of discretizing time has been studied in~\cite{Boland.2019}. Before this summary article, approaches to reducing or avoiding time-expanded networks have been developed in~\cite{Hosseininasab.2015} and~\cite{Boland.2017}. While~\cite{Boland.2017} uses an iterative algorithm using partially time-expanded networks,~\cite{Hosseininasab.2015} models timed events as variables and solves a service network design model without using a time-expanded network at all.
	
During the literature review, it became apparent that there still exists a gap regarding optimization models for timetable-dependent railway network design problems, despite numerous publications in this and neighbouring fields of research. We attempt to reduce this gap by proposing a model for network design models based on strategic timetables. The proposed model uses a multigraph formulation inspired by the one used in~\cite{Sponemann.2013} and avoids the use of time-expanded networks by using timing variables as seen in~\cite{Hosseininasab.2015}. 
	
\section{Methodology}\label{sec:Methodology}
	
In order to build and solve a railway network design problem based on strategic timetables, we need to decide how to model key aspects of railway network design. The following chapters describe three main subjects of input. Section~\ref{sec:ModellingInfrastructure} starts with the modeling of railway infrastructure and capacity, section~\ref{sec:ModellingTimetables} describes how input timetables are modeled and section~\ref{sec:ModellingUncertainty} focuses on uncertainty in the input timetables. After the description of the input data, we present and explain the formulation for the deterministic model (section~\ref{sec:DeterministicModel}) and its robust extension (section~\ref{sec:RobustModel}). Finally, we introduce preprocessing measures to reduce the model's complexity in chapter~\ref{sec:Solving}.
	
\subsection{Modelling Railway Infrastructure}\label{sec:ModellingInfrastructure}
	
We model railway networks as bidirectional multigraphs $G(N,E)$, where nodes $n \in N$ represent railway stations or junctions and arcs $(i,j,tr) \in E$ represent individual tracks on a railway line section $(i,j) \in S$. Using a multigraph allows to include single-, double- or multi-track sections modelled by one parallel arc for each track. Nodes are attributed with a maximal stopping time, which enables the prohibition of stops in nodes which do not feature platforms or overtaking tracks. They also feature individual headway times for single-track travel, which describe the minimal time between the arrival of one train and the departure of a second train on the same track. Headway times of two trains following each other are attributed to the arcs and described in the following paragraph. Since railway stations can vary much in their topology and not all incoming and outgoing lines are necessarily connected, we also consider node links, which represent curves in a railway node connecting incoming and outgoing railway lines. If such a node link is included in the solution, it means that it is possible to travel from nodes $a$ via $i$ to $b$ without changing the direction of travel in $i$. A link is defined by the tuple $(i,a,b)$ with the set $L$ containing all possible links in the network. Node capacities, defined as the number of trains travelling through or dwelling at a node simultaneously, are not limited. This decision is motivated by the difficulties which occur while quantifying a node's capacity on a macroscopic level. The evaluation of a station's capacity and its utilization is heavily dependent on the actual microscopic track layout and the microscopic routing and timing of trains. Since we only consider nodes on a macroscopic level without taking the track layout into account, we are not able to evaluate the station capacities. Please refer to~\cite{Weik.2020} for recent approaches to evaluate the capacity of a railway station while using the framework provided by the UIC Code 406 Leaflet on Railway Capacity~\cite{UIC406}.
	
The arcs $(i,j,tr)$ are identified by the nodes $i$ and $j$, which they connect, and the track number $tr$. Each arc is bidirectional, so trains can use it in both directions. They feature individual and train-type-dependent travel times and train-sequence-dependent headway times. These headway times account for the different travel times and assure that the paths of two trains following each other do not interfere. The headway times are used for the capacity estimation. Please refer to~\cite{UIC406} for an introduction to railway capacity analysis. The model features three different ways to build networks that provide sufficient capacity:
\begin{itemize}
	\item include additional parallel tracks
	\item reduce the given travel times
	\item reduce the given headway times
\end{itemize}
The reduction of travel times covers measures such as the use of different rolling stock with better acceleration, higher velocity or tilting technology or changes to the route like the elimination of tight curves. A reduction of headway times can be achieved by modernizing the signalling system or by changing the amount and position of block signals.

 The routing of trains through the network is path-based. During a preprocessing step, suitable paths are identified for each train in the timetable while respecting travel time bounds and via-nodes. The optimization model chooses exactly one of those paths for each trains. The generation of paths and their impact on the model is further detailed in Section~\ref{sec:Paths}.
	
\subsection{Modelling Input Timetables}\label{sec:ModellingTimetables}
	
The unique speciality of our timetable-based network design problem is that its demand is given by an operational concept derived from an input timetable and not just by train numbers or passenger flows. The input timetable has been created by a strategic timetabling step beforehand and is assumed to be given as input for our model. It includes a macroscopic routing and travel times for all trains. To allow for some optimization of the network, the input timetable has to be relaxed. Therefore only key properties of the strategic timetables are used. The operational concept which becomes the input data for the network design consists of the following elements:
\begin{itemize}
	\item a list of all trains, each one attributed with
	\begin{itemize}
		\item start and destination nodes
		\item optional via-nodes
		\item time bounds given as earliest departure time and latest arrival time
		\item a train type
	\end{itemize}
	\item a set of timing relationships between two trains, which can be
	\begin{itemize}
		\item arrival- or departure-based frequency connections (the arrivals or the departures of two trains of the same line need to be separated by a frequency-time) or
		\item transfer connections (the difference between the arrival of one train and the departure of another has to be within a certain time interval)
	\end{itemize}
\end{itemize}
	
The timing relationships are denoted by triples $(k_1,k_2,i)$ containing the two interacting trains $k_1$ and $k_2$ as well as the node $i$ where the connection is taking place. They each feature a time interval given by $t_{k_1,k_2,i}^{min}$ and $t_{k_1,k_2,i}^{max}$.

\subsection{Modelling Uncertainty}\label{sec:ModellingUncertainty}
	
Both strategic timetabling and network design are planning steps done many years before the actual operation. During the time needed to plan and build infrastructure expansions, political circumstances and demand prognoses can change, which could lead to changes to the strategic timetable. This makes the strategic timetable and the timetable-based network design subject to uncertainty. If a railway network is designed based on one single strategic timetable, an update to this strategic timetable might require adaptions of the network, which may become very expensive. To avoid these additional costs, the probability of these changes should be reduced, and this uncertainty in the strategic timetable should be addressed. Therefore, two different modelling concepts are introduced:
	
\begin{itemize}
	\item Timetable families consisting of discrete timetable scenarios
	\item flexible freight train demand within one scenario, modelled by optional trains
\end{itemize}
	
A timetable family contains several discrete timetable scenarios which are used to model slightly different operational concepts. Each scenario features its list of trains as well as individual timing relationships. Therefore, scenarios might vary in the number of trains, their routes, their train attributes and their timing connections. The robust extension of the model allows the calculation of networks for timetable families. Both fully robust solutions, which cover the whole timetable family, and light robust solutions, where only a defined share of the timetable family is covered, are possible. The grade of robustness is controlled by the scenario share parameter $s_s$ in constraint~(\ref{eq:con:rob:scenario}). This constraints controlles how many percent of the given scenarios have to be covered by the solution network at least. The requested share is given by parameter $s_s$. It is assumed that one scenario is active at any time. Therefore, conflicts may occur between trains of different scenarios.
	
In addition to the timetable families with discrete scenarios, the model also features optional trains, which are used to model uncertain freight train demand. In the deterministic case, all trains included in the list of trains of one scenario are mandatory, so the resulting network has to provide sufficient capacity to run all trains while fulfilling their timing restrictions. In the robust case, this list of mandatory trains can be expanded by a set of optional trains. These optional trains may be included in the solution timetable, but a non-inclusion won't lead to an infeasibility of the model. It is possible to specify a penalty for the non-inclusion of an optional train into the solution. Small penalties will assure, that trains are activated if they fit into the timetable without the need for additional infrastructure. When the values for the penalties are increased, it is possible to create a trade-off between the activation of additional trains and the costs for necessary additional capacity measures. The effect of different ratios between penalties and building costs will be analyzed within a sensitivity analysis, which is going to be part of future research.
	
Besides the penalties and the resulting trade-off, optional trains can also be handled in two other ways. It is possible to demand that a certain number of optional trains has to be included in the solution, while the model chooses which trains can be added at the lowest costs. This is enabled by constraint (\ref{eq:con:rob:opttraincount}). Additionally, the user can also specify certain optional trains to become mandatory. During the optimization, these trains are treated as if they were mandatory in the first place.

\subsection{Deterministic Model}\label{sec:DeterministicModel}
	
The deterministic optimization model features the following binary decision variables ($var \in \left\{0,1\right\}$):

	\begin{align*}
		l_{i,a,b}	&\text{ ... indicates, whether link $(i,a,b)$ is included in the network or not} \\
			p_{k,p}		&\text{ ... indicates, whether path $p$ is used by train $k$ or not} \\
			x_{k,i,j,tr}&\text{ ... indicates, whether train $k$ travels on edge $(i,j,tr)$ or not} \\
			y_{i,j,tr}	&\text{ ... indicates, whether edge $(i,j,tr)$ is included in the network or not} \\
			z_{i,j,tr,k_{1},k_{2}}^{h_c}	&\text{ ... indicates, whether train $k_1$ is running before train $k_2$ on edge $(i,j,tr)$ or not,} \\
			&\quad \text{with $k_1$ crossing $k_2$} \\
			z_{i,j,tr,k_{1},k_{2}}^{h_f}	&\text{ ... indicates, whether train $k_1$ is running before train $k_2$ on edge $(i,j,tr)$ or not,} \\
			&\quad \text{with $k_2$ following $k_1$}
		\end{align*}
	
	Additionally, the model uses the following timing variables, which can be either integers or continuous. All computational results presented later in this paper were calculated with integer timing variables because tests showed that this is the faster choice. Through the use of timing variables, we don't have to use time-expanded networks which greatly reduces the model's size and complexity.
	
	\begin{align*}
			a_{k,i,j,tr}	&\text{ ... arrival of train $k$ when travelling on edge $(i,j,tr)$} \\
			d_{k,i,j,tr}	&\text{ ... departure of train $k$ when travelling on edge $(i,j,tr)$} \\
			r_{i,j}^{mht}	&\text{ ... reduction of minimum headway times on section $(i,j)$, valid for all edges,}\\
			&\text{\quad limited to one minute per ten~km of section length} \\
			r_{i,j}^{time}	&\text{ ... reduction of travel times on section $(i,j)$, valid for all edges, limited so that }\\
			&\text{\quad a defined lower bound for the resulting headway time (e.g. 2 min) cannot be violated	}
		\end{align*}
  
		Besides the variables, several parameters are used by the optimization model:
		\begin{align*}
			a_{k}			&\text{ ... latest possible arrival time of train $k$ at its destination node} \\
			c_{i,j}^{mht}	&\text{ ... specific cost for reducing the headway times on section $(i,j)$ by one minute} \\
			c_{i,j}^{time}	&\text{ ... specific cost for reducing the travel times on section $(i,j)$ by one minute} \\
			d_{k}			&\text{ ... earliest possible departure time for train $k$ at its origin node} \\
			\delta_{k,p}^{i,a,b}&\text{ ... binary parameter which indicates, whether link $(i,a,b)$ is part of path $(k,p)$ or not} \\
			\delta_{k,p}^{i,j}	&\text{ ... binary parameter which indicates, whether section $(i,j)$ is part of path $(k,p)$ or not} \\
			f_{i,a,b}^{link}	&\text{ ... fixed cost for building link $(i,a,b)$} \\
			f_{i,j,tr}^{edge}	&\text{ ... fixed cost for building edge $(i,j,tr)$} \\
			h_{i,j,k_{1},k_{2}}	&\text{ ... minimum headway time between trains $k_1$ and $k_2$ on edge $(i,j,tr)$} \\
			M				&\text{ ... a sufficiently large number (Big-M)} \\
			n_{k}			&\text{ ... number of trains $k$ in the test case} \\
			t_{j}^{cross}	&\text{ ... crossing time in node $j$} \\
			t_{k_{1},k_{2},i}^{con,min}	&\text{ ... lower bound for transfer time between trains $k_1$ and $k_2$ at node $i$} \\
			t_{k_{1},k_{2},i}^{con,max}&\text{ ... upper bound for transfer time between trains $k_1$ and $k_2$ at node $i$} \\
			t_{k_{1},k_{2},i}^{freq,min} &\text{ ... lower bound for time interval between trains $k_1$ and $k_2$ at node $i$} \\
			t_{k_{1},k_{2},i}^{freq,max} &\text{ ... upper bound for time interval between trains $k_1$ and $k_2$ at node $i$} \\
			t_{k,i,j}		&\text{ ... travel time of train $k$ on section $(i,j)$} \\
			t_{maxStop,i}	&\text{ ... maximum stop time in node $i$}
		\end{align*}
 
	The following sets are used for the generation of paths and variables. They are created during the preprocessing, which is described in section~\ref{sec:Solving}
	
		\begin{align*}
			A 	&\text{ ... Set of train-train-node tuples $(k_1,k_2,i)$ for arrival-based frequencies} \\
			C_c &\text{ ... Set of train-train-edge tuples $(k_1,k_2,i,j,tr)$ for crossing headways with conflicts} \\
			D 	&\text{ ... Set of train-train-node tuples $(k_1,k_2,i)$ for departure-based frequencies} \\
			E_a &\text{ ... Set of relevant edges $(i,j,tr)$} \\
			H_c &\text{ ... Set of train-train-edge tuples $(k_1,k_2,i,j,tr)$ for crossing headway cases} \\
			H_f &\text{ ... Set of train-train-edge tuples $(k_1,k_2,i,j,tr)$ for following headway cases} \\
			K 	&\text{ ... Set of trains $k$} \\
			L_a &\text{ ... Set of relevant node links $(i,a,b)$} \\
			N_a &\text{ ... Set of relevant nodes $i$} \\
			O_c &\text{ ... Set of train-train-edge tuples $(k_1,k_2,i,j,tr)$ for crossing headways with fixed order}\\
			O_f &\text{ ... Set of train-train-edge tuples $(k_1,k_2,i,j,tr)$ for following headways with fixed order}\\
			P_c	&\text{ ... Set of train-train-edge tuples $(k_1,k_2,i,j,tr)$ for crossing headway pairs} \\
			P_f &\text{ ... Set of train-train-edge tuples $(k_1,k_2,i,j,tr)$ for following headway pairs} \\
			S_a &\text{ ... Set of relevant line sections $(i,j)$}	\\
			T 	&\text{ ... Set of train-train-node tuples $(k_1,k_2,i)$ for transfers} \\
			TR	&\text{ ... Set of tracks $tr$} \\
			W_k	&\text{ ... Set of paths $p$ for train $k$} \\
			X 	&\text{ ... Set of train-edge tuples $(k,i,j,tr)$}
		\end{align*}
 
	The objective function~(\ref{eq:obj:det}) minimizes infrastructure costs. These include fixed building costs for arcs and links as well as specific costs for each unit of travel or headway time reduction. By setting the building costs to zero, we can model existing infrastructure and shift the focus of the model from network design to network expansion.
	
	\begin{align}
		min \sum\limits_{(i,j,tr)\in E_a} f_{i,j,tr}^{edge}y_{i,j,tr} + \sum\limits_{(i,a,b)\in L_a} f_{i,a,b}^{link}l_{i,a,b} + \sum\limits_{(i,j)\in S_a} c_{i,j}^{time}r_{i,j}^{time} + \sum\limits_{(i,j)\in S_a} c_{i,j}^{mht}r_{i,j}^{mht} \label{eq:obj:det}
	\end{align}
	
	subject to
	
	\begin{align}
		\sum\limits_{(k,p) \in W_k} p_{k,p} &= 1  &&\forall k \in K \label{eq:con:det:paths}\\
		\sum\limits_{(k,p) \in W_k} \delta_{k,p}^{i,a,b}p_{k,p} &\le n_{k}l_{i,a,b} &&\forall (i,a,b) \in L_a
		\label{eq:con:path_link}\\
		\sum\limits_{(k,p) \in W_k} \delta_{k,p}^{i,j}p_{k,p} &= \sum\limits_{tr\in TR}x_{k,i,j,tr} &&\forall k \in K, \forall (i,j) \in S_a \label{eq:con:path_x}\\
		\sum\limits_{k \in K} x_{k,i,j,tr} + x_{k,j,i,tr} &\le n_{k}y_{i,j,tr} &&\forall (i,j,tr) \in E_a \label{eq:con:x_y}\\
		y_{i,j,tr} &\le y_{i,j,tr-1} &&\forall (i,j,tr) \in E_a,\text{ if $tr = 2$} \label{eq:con:track_seq_A} \\
		y_{i,j,tr} &\le y_{i,j,2} &&\forall (i,j,tr) \in E_a,\text{ if $tr > 2$} \label{eq:con:track_seq_B} \\
		d_{k,i,j,tr} + t_{k,i,j} - r_{i,j}^{time} - M\cdot (1-x_{k,i,j,tr}) &= a_{k,i,j,tr} &&\forall (k,i,j,tr) \in X\label{eq:con:traveltime} \\
		d_{k,i,j,tr} + a_{k,i,j,tr} &\le M\cdot x_{k,i,j,tr} &&\forall (k,i,j,tr) \in X \label{eq:con:zerotimes} \\
		\sum\limits_{j\in N_a, tr\in TR} d_{k,o_{k},j,tr} &\ge d_k &&\forall k \in K \label{eq:con:deptimes} \\
		\sum\limits_{i\in N_a, tr\in TR} a_{k,i,e_k,tr} &\le a_k &&\forall k \in K \label{eq:con:arrtimes} \\
		\sum\limits_{l\in N_a,tr\in TR} a_{k,l,i,tr} &\le \sum\limits_{j\in N_a,tr\in TR} d_{k,i,j,tr} &&\forall i \in N_a, \forall k \in K \label{eq:con:nodetiming} \\
		\sum\limits_{j\in N_a,tr\in TR} d_{k,i,j,tr} - \sum\limits_{l\in N_a,tr\in TR} a_{k,l,i,tr} &\le t_{maxStop,i}  &&\forall i \in N_a, \forall k \in K \label{eq:con:maxstoptimes} \\
		r_{i,j}^{time} &= r_{j,i}^{time} &&\forall (i,j) \in S_a \label{eq:con:red_time}\\
		r_{i,j}^{mht} &= r_{j,i}^{mht} &&\forall (i,j) \in S_a \label{eq:con:red_mht} \\
		\sum\limits_{j\in N_a, tr\in TR} a_{k_{2},j,i,tr} - a_{k_{1},j,i,tr} &\ge t_{k_{1},k_{2},i}^{freq,min} &&\forall (k_{1},k_{2},i) \in A \label{eq:con:det:arrfreqmin}\\
		\sum\limits_{j\in N_a, tr\in TR} a_{k_{2},j,i,tr} - a_{k_{1},j,i,tr} &\le t_{k_{1},k_{2},i}^{freq,max} &&\forall (k_{1},k_{2},i) \in A \label{eq:con:det:arrfreqmax}\\
		\sum\limits_{j\in N_a, tr\in TR} d_{k_{2},i,j,tr} - d_{k_{1},i,j,tr} &\ge t_{k_{1},k_{2},i}^{freq,min} &&\forall (k_{1},k_{2},i) \in D \label{eq:con:det:depfreqmin}\\
		\sum\limits_{j\in N_a, tr\in TR} d_{k_{2},i,j,tr} - d_{k_{1},i,j,tr} &\le t_{k_{1},k_{2},i}^{freq,max} &&\forall (k_{1},k_{2},i) \in D \label{eq:con:det:depfreqmax}\\
		\sum\limits_{j\in N_a, tr\in TR} d_{k_{2},i,j,tr} - \sum\limits_{l\in N_a, tr\in TR} a_{k_{1},l,i,tr} &\ge t_{k_{1},k_{2},i}^{con,min} &&\forall (k_{1},k_{2},i) \in T \label{eq:con:det:transfermin}\\
		\sum\limits_{j\in N_a, tr\in TR} d_{k_{2},i,j,tr} - \sum\limits_{l\in N_a, tr\in TR} a_{k_{1},l,i,tr} &\le t_{k_{1},k_{2},i}^{con,max} &&\forall (k_{1},k_{2},i) \in T \label{eq:con:det:transfermax} \\
		z_{i,j,tr,k_{1},k_{2}}^{h_f} &\ge x_{k_1,i,j,tr} + x_{k_2,i,j,tr} - 1 &&\forall (i,j,tr,k_{1},k_{2}) \in O_f \label{eq:con:fixed_seq_one}\\
		x_{k_1,i,j,tr} \cdot x_{k_2,i,j,tr} &= z_{i,j,tr,k_{1},k_{2}}^{h_f}+z_{i,j,tr,k_{2},k_{1}}^{h_f} &&\forall (i,j,tr,k_{1},k_{2}) \in P_f  \label{eq:con:z_mht_1} \\
		h_{i,j,k_{1},k_{2}} - r_{i,j}^{mht} + M\cdot(z_{i,j,tr,k_{1},k_{2}}^{h_f} - 1)  &\le d_{k_{2},i,j,tr}-d_{k_{1},i,j,tr} &&\forall (i,j,tr,k_{1},k_{2}) \in H_f\label{eq:con:mht_1} \\
		x_{k_1,i,j,tr} + x_{k_2,i,j,tr} &\le 1 &&\forall (i,j,tr,k1,k2) \in C_c \label{eq:con:conflict}\\
		z_{i,j,tr,k_{1},k_{2}}^{h_c} &\ge x_{k_1,i,j,tr} + x_{k_2,i,j,tr} - 1 &&\forall (i,j,tr,k_{1},k_{2}) \in O_c \label{eq:con:fixed_seq_two}\\
		x_{k_1,i,j,tr} \cdot x_{k_2,j,i,tr} &= z_{i,j,tr,k_{1},k_{2}}^{h_c}+z_{i,j,tr,k_{2},k_{1}}^{h_c} &&\forall (i,j,tr,k_{1},k_{2}) \in P_c  \label{eq:con:z_mht_2} \\
		CT_{j} + M\times(z_{i,j,tr,k_{1},k_{2}}^{h_c} - 1) &\le d_{k_{2},j,i,tr}-a_{k_{1},i,j,tr}&& \forall (i,j,tr,k_{1},k_{2}) \in H_c\label{eq:con:mht_2}
	\end{align}
		
	The constraints can be grouped into four blocks: constraints (\ref{eq:con:det:paths}) to (\ref{eq:con:track_seq_B}) handle the choice of paths, tracks and node links, (\ref{eq:con:traveltime}) to (\ref{eq:con:red_mht}) deal with travel times and reduction variables, (\ref{eq:con:det:arrfreqmin}) to (\ref{eq:con:det:transfermax}) assure, that all specified connections are respected and (\ref{eq:con:fixed_seq_one}) to (\ref{eq:con:mht_2}) handle the minimal headway times.
	
	For every train $k \in K$ of the operational concept, exactly one path $(k,p)$ has to be chosen (\ref{eq:con:det:paths}). Paths are created beforehand on a node level and are further described in section~\ref{sec:Paths}. The node links $(i,a,b)$ are activated by constraint~(\ref{eq:con:path_link}), if a link belongs to a certain path ($\delta_{k,p}^{i,a,b} = 1$) and the path is activated ($p_{k,p} = 1$). If a certain section $(i,j)$ is part of an active path ($\delta_{k,p}^{i,j} = 1$), one of the available $x$-variables need to be activated (constraint (\ref{eq:con:path_x})), which leads to the activation of the corresponding arc $(i,j,tr)$ by constraint~(\ref{eq:con:x_y}). The x-variables refine the routing by assigning each train to a specific arc. The track-level assignment is important to calculate minimal headway times and estimate capacity usage correctly. Tracks are activated in a certain order. Track two can only be activated after track one, but tracks three and four can be individually activated after track two. Following the rules defined in section~\ref{sec:SymmetryBreaking}, tracks three and four are unidirectional and can therefore be activated if additional capacity is needed for their specific direction.
	
	Besides the network design aspect, the optimization model calculates a feasible macroscopic timetable. Arrival and departure times are assigned to every train-arc combination if the corresponding $x$-variable is active. All of the timing-related constraints contain sums over the arrival or departure variables of all incoming or outgoing arcs for a certain train and a certain node. This is necessary since it is not known, which arc a train will use to get to or leave a node. Since all not-used timing variables are fixed to zero, it can be assumed that only one of the variables within such a sum has a non-zero value.
	
	For each section, a train-type dependant travel time $t_{k,i,j}$ has to be respected (\ref{eq:con:traveltime}). It may be reduced if the reduction variable $r_{i,j}^{time}$ has a non-zero value. If a train does not use a certain track, the associated timing variables are fixed to zero by (\ref{eq:con:zerotimes}). In addition to the travel times, the timing of a train has to respect the given time bounds for departure times $d_k$ at the origin node $o_k$ (\ref{eq:con:deptimes}) and for the arrival times $a_k$ at the destination node $e_k$ (\ref{eq:con:arrtimes}). The node timing constraint (\ref{eq:con:nodetiming}) makes sure, that a train does depart from a station after its arrival. At some nodes, it is not allowed to stop due to a lack of platforms or overtaking tracks. The maximal stopping time $t_{maxStop,i}$ is enforced by (\ref{eq:con:maxstoptimes}). The reduction variables for travel and headway times do always have the same value for both travel directions, which is ensured by constraints (\ref{eq:con:red_time}) and (\ref{eq:con:red_mht}). It is assumed that the improvements reached e.g. by the installation of a modern signalling system are valid for all tracks and both directions of a section. The values of the reduction variables affect the travel time constraint (\ref{eq:con:traveltime}) and the headway time constraint (\ref{eq:con:mht_1}) respectively.
	
	The timetable needs to feature not only correct time allocation for every train but also has to fulfil certain relations between pairs of trains. These timing connections come in three different types. Frequency connections link either two departure or two arrival events of two trains of the same line while transfer connections link the arrival of one train to the departure of another train. The transfer or frequency time is given by an interval with a lower and an upper bound $t_{k_1,k_2,i}^{min}$ and $t_{k_1,k_2,i}^{max}$, which are enforced by a pair of constraints: (\ref{eq:con:det:arrfreqmin}) and (\ref{eq:con:det:arrfreqmax}) for arrival-based frequencies, (\ref{eq:con:det:depfreqmin}) and (\ref{eq:con:det:depfreqmax}) for departure-based frequencies and (\ref{eq:con:det:transfermin}) and (\ref{eq:con:det:transfermax}) for transfers.
	
	The last block of constraints deals with minimal headway times. Different headway situations may occur, and their evaluation based on the time bounds and the shortest travel times is further detailed in section~\ref{sec:Headways}. Headways for two trains travelling on the same track in the same direction are handled by constraints (\ref{eq:con:fixed_seq_one}) to (\ref{eq:con:mht_1}), where (\ref{eq:con:z_mht_1}) selects one of two $z$-variables, which define the sequence of trains $k_1$ and $k_2$. In some cases, the sequence is implicitly fixed by the time bounds. For these situations, the $z$-variable is activated by constraint (\ref{eq:con:fixed_seq_one}) if both trains are actually using the track. The corresponding headway time between the departures of the two trains is then enforced by (\ref{eq:con:mht_1}). For two trains using a track in opposite directions, the sequence decision is handled similarly by (\ref{eq:con:z_mht_2}). Fixed sequences are checked by constraint~(\ref{eq:con:fixed_seq_two}). If it is not possible to route two trains on the same track without violating headway times or time bounds, constraint (\ref{eq:con:conflict}) assures that at most one of them is activated. In the case of opposite directions of travel, the enforcement constraint (\ref{eq:con:mht_2}) connects the departure of the second train with the arrival of the first train. Here, the time to be respected depends not on the trains, but the signalling system in the station, resulting in the crossing time $t_{j}^{cross}$ , which can be individual for each node. 
	
		Constraints~(\ref{eq:con:z_mht_1}) and (\ref{eq:con:z_mht_2}) are formulated as quadratic constraints for better readability. In the code, constraints of the type $a = b \cdot c$, where $a$, $b$ and $c$ are binary terms, are replaced by a set of the following inequalities:
		\begin{align*}
			a &\ge b + c - 1 \\
			a &\le b \\
			a &\le c 
		\end{align*}
		This applies to constraints (\ref{eq:con:z_mht_1}), (\ref{eq:con:z_mht_2}) and (\ref{eq:con:rob:paths}).
			
	\subsection{Robust Model}\label{sec:RobustModel}
	
	In order to account for the uncertainties in the strategic timetable, scenarios and optional trains are used, as described in section~\ref{sec:ModellingUncertainty}. The optimization model had to be expanded to incorporate these additional features. Therefore, the robust model features two additional binary decision variables:
	\begin{itemize}
		\item $o_{s}^{szo}$ - decides, whether scenario $s$ is active or not
		\item $o_{k}^{train}$ - decides, whether train $k$ is active or not
	\end{itemize}
	
		The robust model also requires additional parameters and sets:
		\begin{align*}
			n_{k,opt,s}    &\text{ ... number of optional trains $k$ in scenario $s$} \\
			n_{k,opt,demand,s}  &\text{ ... number of demanded active optional trains $k$ in scenario $s$} \\
			n_{k,s}         &\text{ ... number of trains $k$ in scenario $s$} \\ 
			n_{s}           &\text{ ... number of scenarios $s \in Z$} \\
			p^{szo}_{s}	   &\text{ ... penalty for not activating scenario $s$} \\
			p^{train}_{k}  &\text{ ... penalty for not activating train $k$} \\
			s_{s}          &\text{ ... lower bound of the percentage of scenarios which have to be covered} \\
			& \quad \text{by the resulting network} \\
			K_{mand,s}      &\text{ ... set of mandatory trains $k$ in scenario $s$} \\
			K_{opt,chosen,s}&\text{ ... set of optional trains $k$ in scenario $s$ which have been previously selected} \\
			& \quad \text{to be activated} \\
			Z               &\text{ ... set of all scenarios $s$}
		\end{align*}
	
	The objective function is extended with penalty terms:
	
	\begin{align}
		min \sum\limits_{(i,j,tr)\in E_a} f_{i,j,tr}^{edge}y_{i,j,tr} + \sum\limits_{(i,a,b)\in L_a} f_{i,a,b}^{link}l_{i,a,b} + \sum\limits_{(i,j)\in S_a} c_{i,j}^{time}r_{i,j}^{time} \notag \\ + \sum\limits_{(i,j)\in S_a} c_{i,j}^{mht}r_{i,j}^{mht} + \sum\limits_{k \in K}p^{train}_{k}(1-o_{k}^{train}) + \sum\limits_{s\in Z}p^{szo}_{s}(1-o_{s}^{szo}) \label{eq:obj:rob}
	\end{align}
	
	Additional constraints are introduced to handle scenarios and optional trains:
	
	\begin{align}
		\sum\limits_{s\in Z} \frac{1}{n_{s}}o_{s}^{szo} &\ge s_{s} \label{eq:con:rob:scenario}\\
		\sum\limits_{(k,p) \in P_{k}} p_{k,p} &= o_{k}^{train} \cdot o_{s}^{szo}  &&\forall k \in K \label{eq:con:rob:paths}\\
		o_{k}^{train} &= o_{s}^{szo} &&\forall k \in K_{mand,s} \text{ or } K_{opt,chosen,s}, \forall s \in Z \label{eq:con:rob:trains_szo}\\
		n_{k,s} - n_{k,opt,s} + n_{k,opt,demand,s} &= \sum\limits_{k \in K} o_{k}^{train} &&\forall s \in Z \label{eq:con:rob:opttraincount} \\
		\sum\limits_{j\in N_a, tr\in TR} a_{k_{2},j,i,tr} - a_{k_{1},j,i,tr} &\ge t_{k_{1},k_{2},i}^{freq,min} \cdot o_{s}^{szo} &&\forall s\in Z, \forall (k_{1},k_{2},i) \in A_s \label{eq:con:rob:arrfreqmin}\\
		\sum\limits_{j\in N_a, tr\in TR} a_{k_{2},j,i,tr} - a_{k_{1},j,i,tr} &\le t_{k_{1},k_{2},i}^{freq,max} \cdot o_{s}^{szo} &&\forall s\in Z, \forall (k_{1},k_{2},i) \in A_s \label{eq:con:rob:arrfreqmax}\\
		\sum\limits_{j\in N_a, tr\in TR} d_{k_{2},i,j,tr} - d_{k_{1},i,j,tr} &\ge t_{k_{1},k_{2},i}^{freq,min} \cdot o_{s}^{szo} &&\forall s \in Z, \forall (k_{1},k_{2},i) \in D_s \label{eq:con:rob:depfreqmin}\\
		\sum\limits_{j\in N_a, tr\in TR} d_{k_{2},i,j,tr} - d_{k_{1},i,j,tr} &\le t_{k_{1},k_{2},i}^{freq,max} \cdot o_{s}^{szo} &&\forall s \in Z, \forall (k_{1},k_{2},i) \in D_s \label{eq:con:rob:depfreqmax}\\
		\sum\limits_{j\in N_a, tr\in TR} d_{k_{2},i,j,tr} - \sum\limits_{l\in N_a, tr\in TR} a_{k_{1},l,i,tr} &\ge t_{k_{1},k_{2},i}^{con,min} \cdot o_{s}^{szo} &&\forall s \in Z, \forall (k_{1},k_{2},i) \in T_s\label{eq:con:rob:transfermin}\\
		\sum\limits_{j\in N_a, tr\in TR} d_{k_{2},i,j,tr} - \sum\limits_{l\in N_a, tr\in TR} a_{k_{1},l,i,tr} &\le t_{k_{1},k_{2},i}^{con,max} \cdot o_{s}^{szo} &&\forall s \in Z, \forall (k_{1},k_{2},i) \in T_s
		\label{eq:con:rob:transfermax}
	\end{align}	
	
	The additional constraints (\ref{eq:con:rob:scenario}) to (\ref{eq:con:rob:opttraincount}) handle multiple scenarios as well as optional trains. The important scenario coverage share $s_s$ is enforced by (\ref{eq:con:rob:scenario}), making sure that the percentage of activated scenarios is larger than the coverage share. Constraint (\ref{eq:con:rob:paths}) replaces (\ref{eq:con:det:paths}), because in the robust counterpart, a paths needs to be chosen only if both the scenario and the train are active ($o_k^{train} = o_s^{szo} = 1$). Otherwise, the train does not run, and no path is selected. For the set of mandatory trains within scenario $s$, $K_{mand,s}$, and an optional set of previously chosen optional trains $K_{opt,chosen,s}$, which become mandatory, the $o_k^{train}$-variable has the same value as the correspondent $o_k^{szo}$-variable so that if the scenario is active, the train is active as well. Constraint (\ref{eq:con:rob:opttraincount}) enables the possibility to demand a certain number of optional trains $n_{K,opt,demand,s}$ to be activated while leaving the choice of trains to the solver.
	
	Since the timing relationships can vary among scenarios, they need to be activated and deactivated depending on the scenario status. Therefore, constraints (\ref{eq:con:rob:arrfreqmin}) to (\ref{eq:con:rob:transfermax}) are introduced, which replace (\ref{eq:con:det:arrfreqmin}) to (\ref{eq:con:det:transfermax}). If the scenario is not active, the timing variables are fixed to zero. 
	
	\subsection{Preprocessing measures}\label{sec:Solving}
	
	Early tests indicated that the model is very hard to solve, and optimization times are large. This is primarily caused by a large number of potential train-arc combinations and the consideration of different train succession cases to ensure minimal headway times. The model also suffers from symmetry issues induced by the free choice of tracks for the trains on a section. Several measures were introduced to address these issues and reduce the model's complexity. In section~\ref{sec:Paths}, it is described how paths are used to identify relevant train-arc combinations. Section~\ref{sec:Headways} describes how the number of headway constraints can be reduced by adding constraints only for train pairs where conflicts cannot be resolved beforehand by evaluating the time bounds. Section~\ref{sec:SymmetryBreaking} introduces track-choice rules to deal with symmetry issues and to reduce the number of solutions with the same objective value.
	
	\subsubsection{Introducing paths}\label{sec:Paths}
	
	Initially, the model was formulated purely arc-based, creating $x_{k,i,j,tr}$-variables for every train-arc-combination. This also resulted in many potential headway cases, where each one had to be checked by constraints. However, certain train-arc combinations can be excluded beforehand without losing any solution quality. Two cases can be identified where the combination of train $k$ and arc $(i,j,tr)$ does not need to be included in the optimization model:
	
	\begin{itemize}
		\item the arc $(i,j,tr)$ is not part of any path $p$ connecting origin $o_k$ and destination $e_k$ of train $k$
		\item the minimal travel time on all paths from $o_k$ to $e_k$ that use arc $(i,j,tr)$ is larger than the longest possible travel time for train $k$, leading to the dismissal of the paths
	\end{itemize}
	
	In order to use these cases to reduce the amount of $x$-variables, it is necessary to calculate the available paths for each train beforehand. Therefore, a preprocessing is introduced, which contains two steps:
	
	\begin{itemize}
		\item calculating paths for all (origin,destination,train type)-triples
		\item creating lists of relevant combinations for arcs, trains and headways
	\end{itemize}
	
	The path evaluation is done on a node-section level, meaning that a path contains only a sequence of visited nodes but no information on which of the parallel arcs on a section is used. This is done to limit the number of paths. For the same reason, paths are not created for each train but triples of origin, destination and train type. This way, all trains of the same type that start and end at the same nodes share the same choice of paths. The paths are created before the optimization using the python library \textit{networkx}~\cite{networkx}. To be included in the list of possible paths for a train, a path candidate has to fulfil two requirements:
	\begin{itemize}
		\item its travel time must not exceed the maximal travel time, which is given by the difference between the latest arrival time and the earliest departure time
		\item all specified via nodes must be visited
	\end{itemize}
	
	The network is then reduced to all arcs and node links included in at least one path, resulting in the revised set of arcs $E_a \subseteq E$ and the set of node links $L_a \subseteq L$. The variables $y_{i,j,tr}$ and $l_{i,a,b}$ are created accordingly. After reducing the network, the set of $x$-variables $X$, which includes an quadruple $(k,i,j,tr)$ for each train $k \in K$ and all arcs $(i,j,tr)$ enabled by the calculated paths, is created. This set controls the creation of the $x_{k,i,j,tr}$-variables as well as constraints~(\ref{eq:con:traveltime}) and (\ref{eq:con:zerotimes}), which are only created for the quadruples in $X$ and not for all train-arc-combinations. The constraints~(\ref{eq:con:det:arrfreqmin}) to (\ref{eq:con:det:transfermax}) are created for each entry in the connection sets $A$, $D$ and $T$, indicated by the triple $(k_1,k_2,i)$ denoting the two trains $k_1$ and $k_2$ and the node $i$ where the connection is taking place.
	
	\subsubsection{Evaluating succession cases and time bounds}\label{sec:Headways}
	
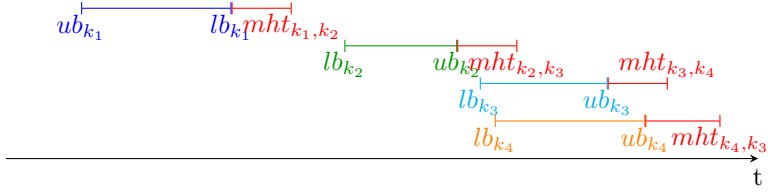
\begin{figure}[H]
		\centering
		\begin{tikzpicture}
			\draw[-stealth] (0,-2.5) -- (10,-2.5);
			\node (t) at (10,-2.75) {t};
			\draw[|-|,blue] (1,-0.5) -- (3,-0.5);
			\draw[|-|,red] (3,-0.5) -- (3.8,-0.5);
			\node[blue] (ub_k1) at (1,-0.75) {$ub_{k_1}$};
			\node[blue] (lb_k1) at (3,-0.75) {$lb_{k_1}$};
			\node[red] (mht_k1k2) at (3.8,-0.75) {$mht_{k_1,k_2}$};
			\draw[|-|,darkgreen] (4.5,-1) -- (6,-1);
			\draw[|-|,red] (6,-1) -- (6.8,-1);
			\node[darkgreen] (ub_k2) at (4.5,-1.25) {$lb_{k_2}$};
			\node[darkgreen] (lb_k2) at (6,-1.25) {$ub_{k_2}$};
			\node[red] (mht_k2k3) at (6.8,-1.25) {$mht_{k_2,k_3}$};
			\draw[|-|,cyan] (6.3,-1.5) -- (8,-1.5);
			\draw[|-|,red] (8,-1.5) -- (8.8,-1.5);
			\node[cyan] (lb_k3) at (6.3,-1.75) {$lb_{k_3}$};
			\node[cyan] (ub_k3) at (8,-1.75) {$ub_{k_3}$};
			\node[red] (mht_k3k4) at (8.8,-1.25) {$mht_{k_3,k_4}$};
			\draw[|-|,orange] (6.5,-2) -- (8.5,-2);
			\draw[|-|,red] (8.5,-2) -- (9.5,-2);
			\node[orange] (lb_k4) at (6.5,-2.25) {$lb_{k_4}$};
			\node[orange] (ub_k4) at (8.5,-2.25) {$ub_{k_4}$};
			\node[red] (mht_k4k3) at (9.5,-2.25) {$mht_{k_4,k_3}$};
		\end{tikzpicture}
		\caption{Comparison of time windows to identify headway cases}
		\label{fig:timewindowsheadway}
	\end{figure}

	Making sure that all headway times are correctly respected has a large impact on the model's size and computation times. To reduce the very large number of train-train-arc combinations, potential headways are evaluated before building the optimization model. The goal is to identify those succession cases, where headway times have to be explicitly ensured by constraints in the optimization model and to discard all cases, where trains cannot use a certain track together or headway times are implicitly assured.
	For the headway time constraints~(\ref{eq:con:z_mht_1}) to (\ref{eq:con:mht_2}), several sets are needed. The sets $P_f$ and $P_c$ include quintuples $(i,j,tr,k_1,k_2)$ and contain all pairs of trains $k_1$ and $k_2$ which might use the same track $(i,j,tr)$ following ($P_f$) or crossing ($P_c$) each other. The sequence of the two trains is not relevant for these sets but gets considered in sets $H_f$ and $H_c$. Based on the time bounds $d_{k_1}$ and $d_{k_2}$ for the earliest departure time at the start node and $a_{k_1}$ and $a_{k_2}$ for the latest arrival time at the destination node and the minimal travel times between origin, destination and the start node of the edge in question, lower and upper bounds $lb_{k_1}$, $ub_{k_1}$, $lb_{k_2}$ and $ub_{k_2}$ for the departure of both trains can be calculated. The different possible intersections of these time windows define, whether headway time constraints need to be created or not. Several possible cases are depicted in figure~\ref{fig:timewindowsheadway}.
	If two time windows are disjoint and the gap between the time windows is larger than the headway time, the headway time is implicitly respected and does not need to be secured by the optimization model. This is the case for trains $k_1$ and $k_2$ in figure~\ref{fig:timewindowsheadway}. Between trains $k_2$ and $k_3$, the sequence is fixed, because $ub_{k_2} < lb_{k_1}$. But since $lb_{k_3} < ub_{k_2} + mht_{k_2,k_3}$, the headway time needs to be assured by a constraint. Therefore, the quintuple$(i,j,tr,k_2,k_3)$ is added to set $O_f$ or $O_c$ and $H_f$ or $H_c$ and a constraint pair (\ref{eq:con:fixed_seq_one}), (\ref{eq:con:mht_1}) or (\ref{eq:con:fixed_seq_two}), (\ref{eq:con:mht_2}) is created. The time windows of trains $k_3$ and $k_4$ overlap so much that both sequences could be feasible. In this case, constraint pairs (\ref{eq:con:z_mht_1}), (\ref{eq:con:mht_1}) or (\ref{eq:con:z_mht_2}), (\ref{eq:con:mht_2}) are needed to determine the sequence and the actual departure times. 
	For a special case of two trains travelling in opposite directions, conflicts can be detected by evaluating the time windows. This happens if the sequence is fixed to $k_1$ before $k_2$ at node $i$ on a section $(i,j)$ and to $k_2$ before $k_1$ at node $j$. This implies that the two trains have to cross each other on this section, which would lead to a collision if both trains are using the same track. Conflicted pairs are added to set $C_c$ and these conflicts are avoided by constraint~(\ref{eq:con:conflict}).
	
	\subsubsection{Breaking symmetry by introducing track choice rules}\label{sec:SymmetryBreaking}
	
	Early versions of the optimization model used all train-arc combinations for sections available to a train. This implicated that trains can freely choose between all available tracks on a section. Combining this with the potentially large number of sections travelled on during a train's trip led to many solutions with equal objective values, which all needed to be evaluated to prove optimality. The consequence was that the optimization took much time, often without improving the objective value. The existence of large numbers of equal-valued solutions, originating from the ability to permute variables without changing the structure of the problem, in an (integer) linear program is also known as symmetry~\cite{Margot.2010}. It is possible to reduce this symmetry by introducing special rules or constraints. In this case, rules limiting the choice of tracks for each train were introduced. Consider a section $A - B$ and a maximum number of tracks of two. This section is used by various trains travelling from A to B or from B to A. If only a few trains use the section, one track may provide enough capacity for all trains. In this case, only the building costs for one track must be paid, and all trains use this track.
	
\begin{figure}[H]
		\centering
		\begin{tikzpicture}
			\node[draw=black,minimum height = 1cm] (A1) at (0,0) {A};
			\node[draw=black,minimum height = 1cm] (B1) at (3,0) {B};
			\draw[stealth-stealth, black, thick] (A1.east) -- (B1.west);
			\draw[-stealth,blue] ([yshift=-0.2cm]A1.east) -- ([yshift=-0.2cm]B1.west);
			\draw[stealth-,red] ([yshift = 0.2cm]A1.east) -- ([yshift=0.2cm]B1.west);
			\node[draw=black,minimum height = 2cm] (A2) at (4,0) {A};
			\node[draw=black,minimum height = 2cm] (B2) at (7,0) {B};
			\draw[stealth-stealth, black, thick] ([yshift=-0.5cm]A2.east) -- ([yshift=-0.5cm]B2.west);
			\draw[-stealth,blue] ([yshift=-0.7cm]A2.east) -- ([yshift=-0.7cm]B2.west);
			\draw[stealth-,red] ([yshift =-0.3cm]A2.east) -- ([yshift=-0.3cm]B2.west);
			\draw[stealth-stealth, black, thick] ([yshift=0.5cm]A2.east) -- ([yshift=0.5cm]B2.west);
			\draw[stealth-,red] ([yshift =0.7cm]A2.east) -- ([yshift=0.7cm]B2.west);
			\node[draw=black,minimum height = 4cm] (A4) at (8,0) {A};
			\node[draw=black,minimum height = 4cm] (B4) at (11,0) {B};
			\draw[stealth-stealth, black, thick] ([yshift=-1.5cm]A4.east) -- ([yshift=-1.5cm]B4.west);
			\draw[-stealth,blue] ([yshift=-1.7cm]A4.east) -- ([yshift=-1.7cm]B4.west);
			\draw[stealth-,red] ([yshift =-1.3cm]A4.east) -- ([yshift=-1.3cm]B4.west);
			\draw[stealth-stealth, black, thick] ([yshift=-0.5cm]A4.east) -- ([yshift=-0.5cm]B4.west);
			\draw[stealth-,red] ([yshift =-0.3cm]A4.east) -- ([yshift=-0.3cm]B4.west);
			\draw[stealth-stealth, black, thick] ([yshift=0.5cm]A4.east) -- ([yshift=0.5cm]B4.west);
			\draw[-stealth,blue] ([yshift=0.3cm]A4.east) -- ([yshift=0.3cm]B4.west);
			\draw[stealth-stealth, black, thick] ([yshift=1.5cm]A4.east) -- ([yshift=1.5cm]B4.west);
			\draw[stealth-,red] ([yshift =1.7cm]A4.east) -- ([yshift=1.7cm]B4.west);
		\end{tikzpicture}
		\caption{Track choices for different track numbers for trains of different directions}
		\label{fig:trackchoices}
	\end{figure}
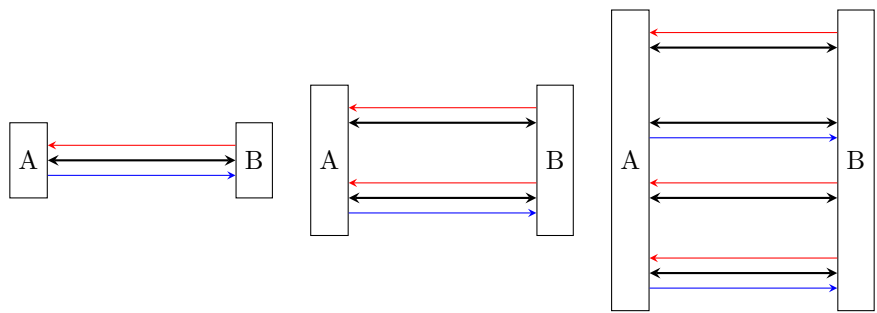

	However, as soon as a second track is needed, all trains have the option to use both tracks and many different solutions are created, limited only by the other trains and the necessary headways between a pair of trains. The objective value is the same for all cases since it is affected only by the building costs for the two tracks. 
	To reduce the number of potential solutions, two rules have been established:
	\begin{itemize}
		\item All trains travelling in alphabetically ascending direction ($A - B$) use track one
		\item All trains travelling in alphabetically descending direction ($B - A$) can choose between tracks one and two
	\end{itemize}
	These rules still allow for creating both single- or double-track sections. However, the amount of symmetric solution is greatly reduced. These rules have proven to be highly effective in reducing the optimization time. They also increase the realism of the obtained solutions since, in most countries, trains on double-track lines are sorted by direction, e.g. in Germany, trains usually use the right track in their direction of travel, while in Switzerland, trains use the left track. The different choices for trains travelling $A - B$ (blue) and $(B - A)$ (red) are illustrated in figure~\ref{fig:trackchoices}. In sections with more than two tracks, the rules are expanded so that trains travelling in ascending order use the odd track numbers, while trains travelling in descending order can either use track one or a track with even number. Constraint~(\ref{eq:con:track_seq_B}) guarantees, that tracks three and four are only built after track two, but they do not depend on each other. More than four tracks are possible, but this case is very rare in praxis and therefore not considered.
	
	\section{Case study and computational results}\label{sec:CaseStudy}
	
	To demonstrate the model, identify performance issues and prove its functionality, a small case study has been conducted. Both the deterministic and the robust optimization models described in section~\ref{sec:DeterministicModel} and \ref{sec:RobustModel} have been implemented in Python 3.8 and solved with Gurobi 10.0.1 All test cases have been solved on a personal laptop featuring an Intel i7-8565U CPU running at 1.80 GHz and 16 GB of RAM. Using the publicly available timetable drafts from the \textit{Deutschlandtakt}~\cite{SMAundPartnerAG.2020.ITF-SO}, twelve input timetables have been created. They consist of up to 716 different trains of up to five train types, travelling through a network of up to 125 nodes and 660 edges. This section provides insight into test cases and computational results for both the deterministic model (section~\ref{sec:CompResDet}) and the robust model (section~\ref{sec:CompResRobust}).
	
	\subsection{Test cases for the deterministic, single-scenario optimization model}
	\label{sec:CompResDet}
	
	The number of different trains per type and timing connections of the test instances are summarized in table~\ref{tab:TestCasesDeterministic}, where the trains are grouped by their train type, either highspeed (FV), intercity (FR), regional express (X), local (N) or freight (G). All the smaller test cases are subset of the largest one \textit{Saxony-716-4h}. Each of these test cases represents one scenario, but they can be expanded to a timetable family for the robust case by creating several scenarios. The test case identifier includes the region, the number of trains and the timespan during which all trains depart from their origin node.
	
	\begin{table}[H]
		\centering
		\caption{Test cases for the deterministic model}
		\label{tab:TestCasesDeterministic}
			\begin{tabular}{llllllll}
				\textbf{Test Case} & \textbf{FV} & \textbf{FR} & \textbf{X} & \textbf{N} & \textbf{G} & \textbf{Connections} & \textbf{Frequencies}\\
				\hline
				Eastern-Saxony-local-22-2h & 0 & 0 & 0 & 22 & 0 & 12 & 32\\
				Eastern-Saxony-local-64-2h & 0 & 0 & 0 & 64 & 0 & 14 & 56\\
				Eastern-Saxony-28-1h & 4 & 4 & 0 & 12 & 8 & 0 & 4 \\
				Eastern-Saxony-40-2h & 6 & 2 & 0 & 22 & 10 & 3 & 24 \\
				Eastern-Saxony-92-2h & 6 & 6 & 0 & 60 & 20 & 17 & 64 \\
				Eastern-Saxony-184-4h & 12 & 12 & 0 & 120 & 40 & 60 & 226 \\
				Saxony-Local-2h & 0 & 0 & 28 & 248 & 0 & 32 & 313 \\
				Saxony-Intercity-2h & 6 & 8 & 0 & 0 & 0 & 4 & 8 \\
				Saxony-Freight-2h & 0 & 0 & 0 & 0 & 68 & 0 & 0 \\
				Saxony-174-1h & 4 & 4 & 15 & 122 & 33 & 10 & 58 \\ 
				Saxony-358-2h & 6 & 8 & 28 & 248 & 68 & 38 & 321 \\ 
				Saxony-716-4h & 12 & 16 & 56 & 496 & 136 & 38 & 831 \\
			\end{tabular}
	\end{table}
	
	All test cases have been solved by Gurobi. The results, including the time needed for the preprocessing, are shown in table~\ref{tab:ComputationalResultsDeterministic}. The preprocessing consists of the path generation, the creation of the various sets of relevant combinations including the headway time window evaluation and the generation of the optimization model. The time limit for the optimization was set to two hours. As mentioned in section~\ref{sec:ModellingInfrastructure}, the model features three options to provide capacity: activating additional parallel edges or reducing travel and headway times. To evaluate their influence, each model has been solved using three different configurations:
	\begin{description}
		\item[Config A] up to four parallel edges, time reductions allowed
		\item[Config B] up to two parallel edges, time reductions allowed
		\item[Config C] up to two parallel edges, time reductions prohibited
	\end{description}
	
	\begin{table}[H]
		\centering
		\caption{Computational results for the deterministic model}
		\label{tab:ComputationalResultsDeterministic}
		\resizebox{\linewidth}{!}{
			\begin{tabular}{lllllllll}
				\textbf{Test Case}       & \textbf{Variables} & \textbf{Constraints}   & \textbf{Building}&\textbf{Edges} &\textbf{Links} & \textbf{Preprocessing} & \textbf{Runtime [s]}   & \textbf{Gap [\%]}  \\
				 &   &  & \textbf{Costs} & & & \textbf{Time [s]}      &     &           \\
				\hline
				Eastern-Saxony-local-22-2h (A) & 1985 & 3431 & 382300 & 14 & 8 & 0.83 & 0.37 & 0.00 \\
				Eastern-Saxony-local-22-2h (B) & 1332 & 2687 & 382300 & 14 & 8 & 0.49 & 0.12 & 0.00 \\
				Eastern-Saxony-local-22-2h (C) & 1098 & 2453 & 382300 & 14 & 8 & 0.60 & 0.07 & 0.00 \\
				Eastern-Saxony-local-64-2h (A) & 6178 & 13870 & 1335300 & 57 & 33 & 2.83 & 78.49 & 0.00 \\
				Eastern-Saxony-local-64-2h (B) & 4158 & 11657 & 1335300 & 57 & 33 & 2.46 & 7.52 & 0.00 \\
				Eastern-Saxony-local-64-2h (C) & 3340 & 10839 & 1335300 & 57 & 33 & 2.41 & 6.11 & 0.00 \\
				Eastern-Saxony-28-1h (A) & 11654 & 24613 & 354000 & 29 & 25 & 4.32 & 99.35 & 0.00 \\
				Eastern-Saxony-28-1h (B) & 7965 & 18257 & 354000 & 29 & 25 & 3.02 & 51.46 & 0.00 \\
				Eastern-Saxony-28-1h (C) & 7216 & 17509 & 354000 & 29 & 25 & 2.42 & 42.51 & 0.00 \\
				Eastern-Saxony-40-2h (A) & 19776 & 38340 & 800000 & 49 & 51 & 57.44 & 7200.00 & 8.91 \\
				Eastern-Saxony-40-2h (B) & 13669 & 28972 & 779800 & 46 & 44 & 41.18 & 7200.00 & 2.46 \\
				Eastern-Saxony-40-2h (C) & 12004 & 27307 & 779800 & 46 & 44 & 36.03 & 4886.55 & 0.00 \\
				Eastern-Saxony-92-2h (A) & 33921 & 73262 & 1779000 & 85 & 64 & 45.46 & 7200.00 & 5.97 \\
				Eastern-Saxony-92-2h (B) & 23434 & 57127 & 1771200 & 86 & 63 & 61.44 & 7200.00 & 5.49 \\
				Eastern-Saxony-92-2h (C) & 20951 & 54644 & 1762400 & 84 & 66 & 47.91 & 7200.00 & 3.53 \\
				Eastern-Saxony-184-4h (A) & 73863 & 163073 & 2302000 & 111 & 84 & 87.60 & 7200.00 & 26.04 \\
				Eastern-Saxony-184-4h (B) & 51220 & 127502 & 1885000 & 94 & 64 & 75.80 & 7200.00 & 9.75 \\
				Eastern-Saxony-184-4h (C) & 46566 & 122848 & - & - & - & 127.00 & 7200.00 & - \\
				Saxony-Local-2h (A) & 69463 & 202078 & - & - & - & 66.1 & 7200.00 & - \\
				Saxony-Local-2h (B) & 48241 & 167117 & 3587900 & 189 & 123 & 47.99 & 7200.00 & 4.76 \\
				Saxony-Local-2h (C) & 43865 & 162741 & 3562000 & 186 & 131 & 44.87 & 7200.00 & 2.31 \\
				Saxony-Intercity-2h (A) & 5447 & 6315 & 588900 & 37 & 31 & 8.54 & 16.43 & 0.00 \\
				Saxony-Intercity-2h (B) & 3705 & 4683 & 588900 & 37 & 31 & 8.54 & 2.22 & 0.00 \\
				Saxony-Intercity-2h (C) & 2806 & 3784 & 588900 & 37 & 31 & 12.70 & 3.49 & 0.00 \\
				Saxony-Freight-2h (A) & 37264 & 75218 & 1118400 & 69 & 60 & 61.84 & 7200.00 & 6.36 \\
				Saxony-Freight-2h (B) & 25393 & 57205 & 1108000 & 65 & 56 & 51.69 & 7200.00 & 2.54 \\
				Saxony-Freight-2h (C) & 22490 & 54302 & 1108000 & 65 & 56 & 57.07 & 7200.00 & 0.12 \\
				Saxony-174-1h (A) & 66515 & 171630 & 3446800 & 191 & 154 & 83.2 & 7200.00 & 10.01 \\
				Saxony-174-1h (B) & 46202 & 137587 & 3254400 & 173 & 146 & 77.98 & 7200.00 & 2.89 \\
                Saxony-174-1h (C) & 42139 & 133524 & 3241200 & 172 & 148 & 79.90 & 7200.00 & 1.30 \\
				Saxony-358-2h (A) & 162800 & 433206 & - & - & - & 177.70 & 7200.00 & - \\
				Saxony-358-2h (B) & 114172 & 344986 & - & - & - & 135.70 & 7200.00 & - \\
				Saxony-358-2h (C) & 106746 & 337560 & - & - & - & 167.90 & 7200.00 & - \\
				Saxony-716-4h (A) & 364305 & 975944 & - & - & - & 423.70 & 7200.00 & - \\
				Saxony-716-4h (B) & 257231 & 775516 & - & - & - & 306.30 & 7200.00 & - \\
				Saxony-716-4h (C) & 243021 & 761306 & - & - & - & 338.00 & 7200.00 & - \\
			\end{tabular}}
	\end{table}
	
	The computational results in table~\ref{tab:ComputationalResultsDeterministic} show, that only very small instances can be solved to optimality within 2~hours. The computation times depend not only on the number of trains in an instance, but also on the train types, the length of train routes, the number of available paths and the flexibility within the time bounds of each train. It gets harder to solve the model and to prove an optimal solution with more train types interacting with each other, with trains travelling longer routes and having more paths to chose from and with larger flexibility within the time bounds. When comparing the different configurations, it becomes apparent that reducing the number of parallel edges and removing the reduction variables leads to smaller models that are a little easier to solve without giving up solution quality. This leads to assumption that the reduction variables are seldom used in an optimal solution. For the two largest instances \textit{Saxony-358-2h} and \textit{Saxony-716-4h} and configuration C of instance \textit{Eastern-Saxony-184-4h}, the solver was unable to find a feasible solution within the given time frame. This indicates that improvements to the solvers performance, for example through the use of decomposition algorithms, are necessary to solve these or even larger instances within reasonable time, even though computation times are not the highest priority in a strategic setting.
	
\begin{figure}[H]
		\centering
		\includegraphics[width=10cm]{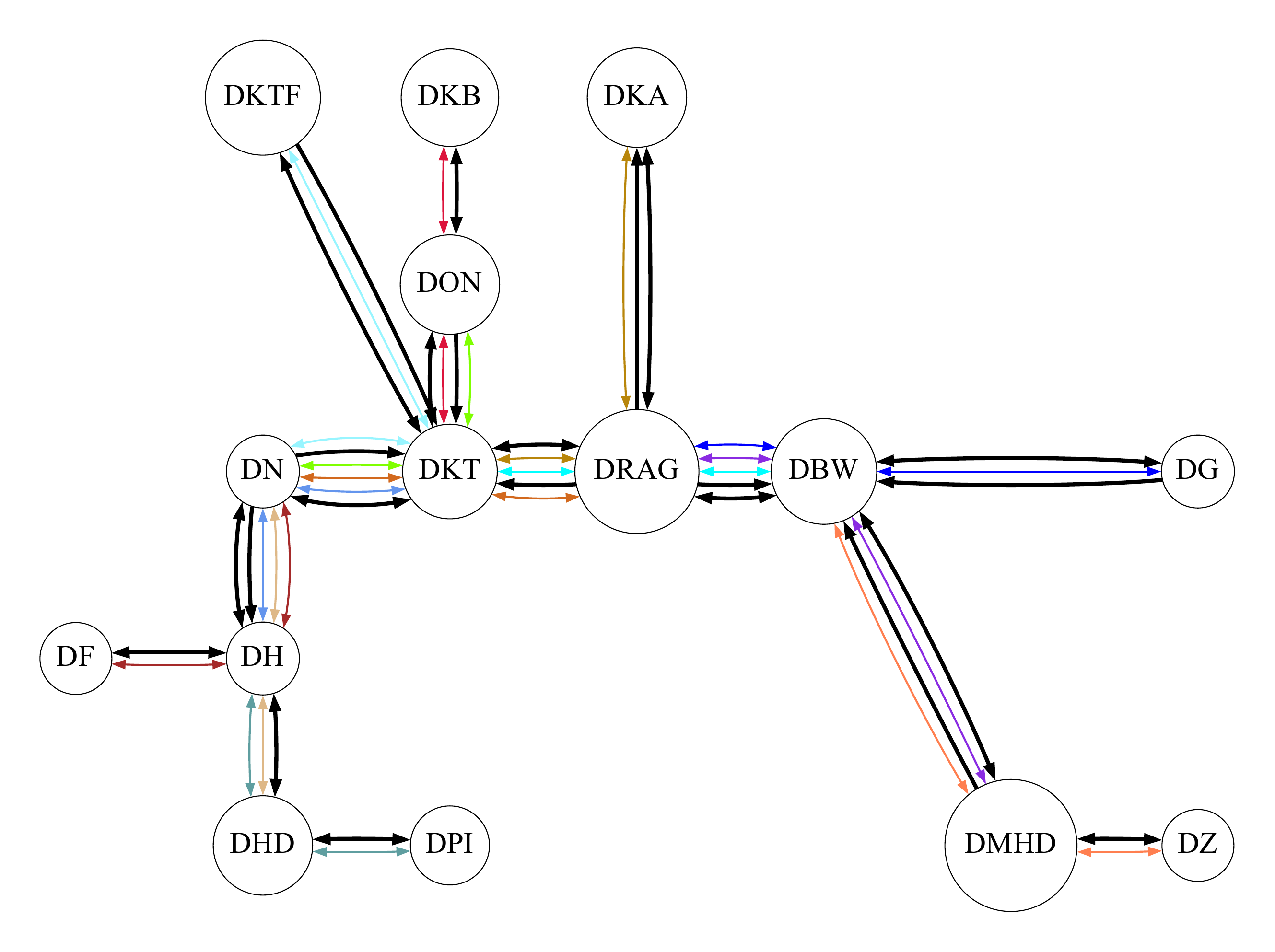}
		\caption{Exemplary solution: infrastructure}
		\label{fig:exNetwork}
	\end{figure}

\begin{figure}[H]
		\centering
		\includegraphics[width=10cm]{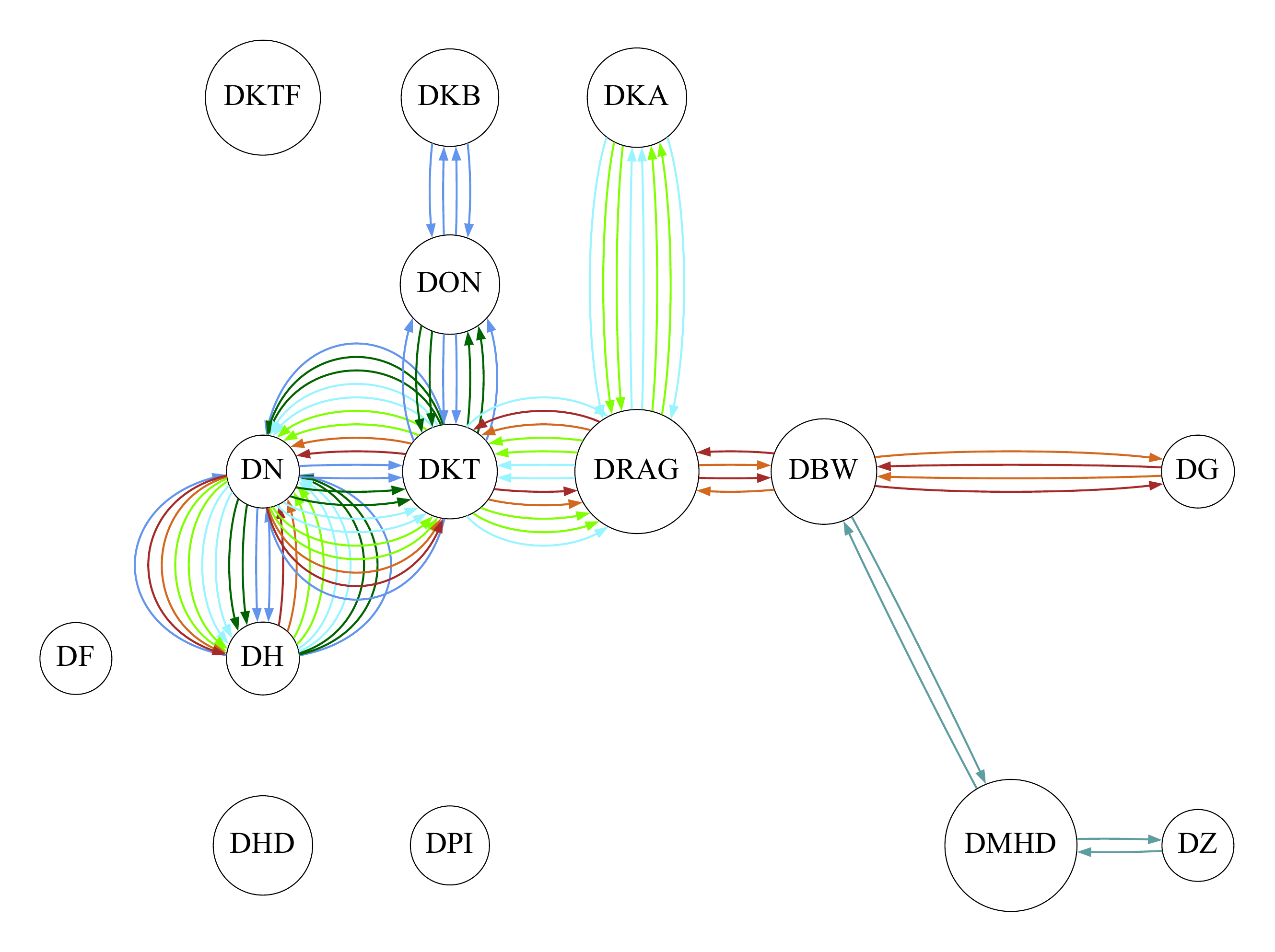}
		\caption{Exemplary solution: train routing}
		\label{fig:exRouting}
	\end{figure}

\begin{figure}[H]
		\centering
		\includegraphics[width=10cm]{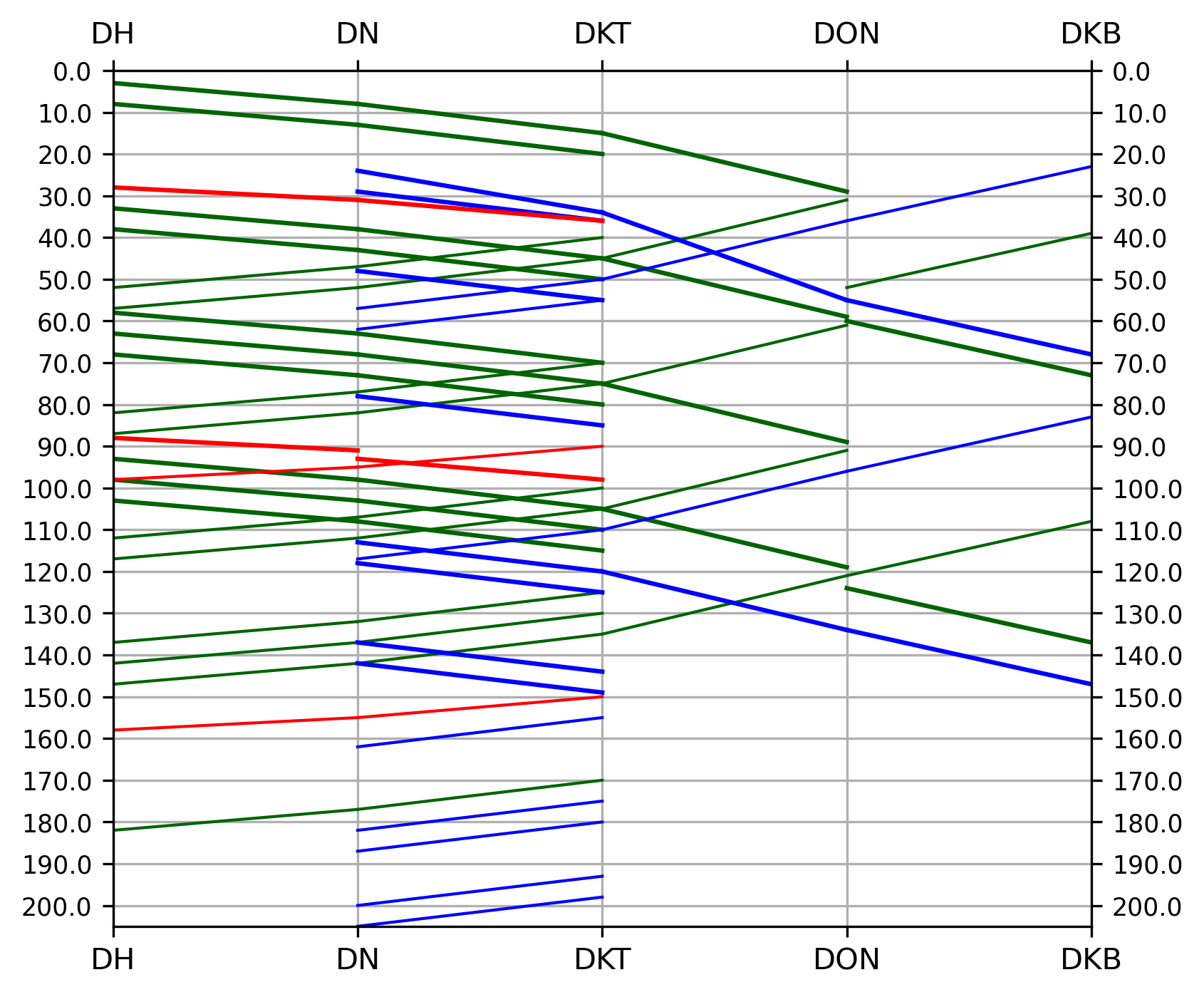}
		\caption{Exemplary solution: macroscopic timetable}
		\label{fig:exTimetable}
	\end{figure}

	A typical solution of the optimization model consists of three components: the network (see figure~\ref{fig:exNetwork}), the routing (see figure~\ref{fig:exRouting}) and the macroscopic timetable (see figure~\ref{fig:exTimetable}). The network contains all nodes, arcs and links necessary to provide the capacity for all demanded trains. In figure~\ref{fig:exNetwork} the arcs are indicated by bold, black arrows while a link $l_{i,a,b}$ is denoted by a pair of coloured arrows from $a$ to $i$ and from $i$ to $b$. The direction of the arrows matches the track-choice rules defined in Section~\ref{sec:SymmetryBreaking}, so each section features one bi-directional arc, while additional arcs are unidirectional. All scenario trains are routed through the network as depicted in figure~\ref{fig:exNetwork}. The trains are coloured by their lines, so trains travelling on the same route are depicted in the same colour. Figure~\ref{fig:exNetwork} shows a solution from a robust test case, which is why the network contains nodes that are not visited by any of the trains in the depicted scenario. However, these are part of another scenario which the network covers. Besides the routing, also the timing is fixed during the optimization. The timing contains arrival and departure times for all nodes visited by the trains and is calculated in a way that fulfils the timing restrictions. The resulting macroscopic timetable can be depicted as a time-space-diagram as illustrated in figure~\ref{fig:exTimetable}. Here, the colours represent different train types: intercity trains are drawn in red, local trains in green, and freight trains in blue.
	
	\subsection{Test cases and computational results for the robust, multi-scenario model} \label{sec:CompResRobust}
	
	For the robust case, two of the smaller test cases \textit{Eastern-Saxony-local-22-2h} and \textit{Eastern-Saxony-40-2h} have been expanded to a timetable family which contains ten different scenarios. They vary in the type, route and amount of included trains and the demanded frequencies and connections. The different scenarios are summarised in table~\ref{tab:ScenarioOverview22} for \textit{Eastern-Saxony-local-22-2h} and in table~\ref{tab:ScenarioOverview40} for \textit{Eastern-Saxony-40-2h}.
	
	\begin{table}[H]
		\centering
		\caption{Scenario overview for the robust test case \textit{Eastern-Saxony-Local-22-2h}}
		\label{tab:ScenarioOverview22}
			\resizebox{\linewidth}{!}{
			\begin{tabular}{llllll}
				\textbf{Scenario} & \textbf{Intercity} & \textbf{Local} & \textbf{Freight} & \textbf{Connections} & \textbf{Frequencies} \\
				& \textbf{Trains} & \textbf{Trains} & \textbf{Trains} & & \\
				\hline
				Eastern-Saxony-local-22-0 & 0 & 22 & 0  & 0  & 0\\
				Eastern-Saxony-local-22-1 & 0 & 22 & 0  & 12 & 32\\
				Eastern-Saxony-local-22-2 & 0 & 22 & 0  & 12 & 32\\
				Eastern-Saxony-local-22-3 & 0 & 22 & 4 (Route A) & 12 & 32\\
				Eastern-Saxony-local-22-4 & 0 & 22 & 4 (Route B) & 12 & 32\\
				Eastern-Saxony-local-22-5 & 0 & 22 & 4 (Route C) & 12 & 32\\
				Eastern-Saxony-local-22-6 & 0 & 22 & 4 (Route D) & 12 & 32\\
				Eastern-Saxony-local-22-7 & 0 & 22 & 16 & 12 & 32\\
				Eastern-Saxony-local-22-8 & 6 & 22 & 16 & 14 & 38\\
				Eastern-Saxony-local-22-9 & 0 & 30 & 0  & 12 & 44\\
			\end{tabular}}
	\end{table}
	
	\begin{table}[H]
		\centering
		\caption{Computational results for the robust test case \textit{Eastern-Saxony-Local-22-2h}}
		\label{tab:ComputationalResultsRobust22}
			\begin{tabular}{lllllll}
				\textbf{Coverage}    & \textbf{Active}        & \textbf{Infrastructure}    & \textbf{Arcs}  & \textbf{Links} & \textbf{Runtime [s]}   & \textbf{Gap [\%]}  \\
				\textbf{Share [\%]}  & \textbf{Scenarios [\%]}& \textbf{Costs}             &  &    &               &           \\
				\hline
				10 & 20 & 393300 & 15 & 8 & 20.60 & 0.00 \\
				20 & 30 & 397100 & 16 & 8 & 14.34 & 0.00 \\
				30 & 30 & 397100 & 16 & 8 & 37.81 & 0.00 \\
				40 & 60 & 400900 & 17 & 9 & 61.60 & 0.00 \\
				50 & 60 & 400900 & 17 & 9 & 62.61 & 0.00 \\
				60 & 60 & 400900 & 17 & 9 & 53.03 & 0.00 \\
				70 & 70 & 401900 & 17 & 10 & 155.46 & 0.00 \\
				80 & 80 & 411700 & 18 & 10 & 112.96 & 0.00 \\
				90 & 90 & 432400 & 19 & 10 & 611.29 & 0.00 \\
				100 & 100 & 458400 & 23 & 13 & 60.09 & 0.00 \\
			\end{tabular}
	\end{table}
	
	With all these scenarios, networks covering different shares of the scenarios have been calculated using the expansion measures defined by configuration B. Therefore, the parameter $s_{s}$ in constraint~(\ref{eq:con:rob:scenario}) has been varied in steps of ten per cent. The main results and computation times can be found in table~\ref{tab:ComputationalResultsRobust22} for \textit{Eastern-Saxony-local-22-2h} and in table~\ref{tab:ScenarioOverview40} for \textit{Eastern-Saxony-40-2h}. Several effects are apparent: first, the infrastructure costs and the amount of included arcs in the network increase with increasing coverage share. This matches expectations since more scenarios to cover also mean more varied demand that has to be satisfied. Second, on some occasions, the scenario coverage share is over-fulfilled because several scenarios can be activated without including additional arcs and causing additional infrastructure costs. Since the robust model with ten scenarios has more variables (24788 instead of 1332 and 92091 instead of 13669) and constraints (64917 compared to 2687 and 223764 compared to 28972) than the deterministic model with only one scenario, it is to be expected that the robust model takes longer to solve. However, it is interesting to see that the computation times vary a lot depending on the value of the coverage share. The model is easier to solve for low coverage shares and for the full robust solution, while high coverage shares lead to longer computation times. For the second test case \textit{Eastern-Saxony-40}, no optimal solution could be calculated within two hours, regardless of the coverage share. The dramatic increase in the computation times compared to the deterministic model illustrates once again, that performance improvements and also heuristic approaches for the robust model are necessary to be able to compute robust solutions for larger test cases.
    Nevertheless, calculating robust, timetable-based solutions combines the advantages of a cost-efficient network, that is optimized towards a given timetable, with the flexibility of being able to choose from one of the feasible timetable scenarios in a later planning stage. This planning method makes sure, that the timetable calculation is not overly restricted by the infrastructure, which might happen if the network would be planned without taking the timetable into account. 
	
	\begin{table}[H]
		\centering
		\caption{Scenario overview for the robust test case \textit{Eastern-Saxony-40-2h}}
		\label{tab:ScenarioOverview40}
			\begin{tabular}{llllll}
				\textbf{Scenario} & \textbf{Intercity} & \textbf{Local} & \textbf{Freight} & \textbf{Connections} & \textbf{Frequencies} \\
				& \textbf{Trains} & \textbf{Trains} & \textbf{Trains} & & \\
				\hline
				Eastern-Saxony-40-0 & 8 & 22 & 10  & 3  & 24\\
				Eastern-Saxony-40-1 & 8 & 22 & 10  & 10 & 24\\
				Eastern-Saxony-40-2 & 8 & 26 & 10  & 3 & 36\\
				Eastern-Saxony-40-3 & 8 & 26 & 10  & 3 & 36\\
				Eastern-Saxony-40-4 & 8 & 30 & 10  & 3 & 48\\
				Eastern-Saxony-40-5 & 8 & 30 & 10  & 15 & 48\\
				Eastern-Saxony-40-6 & 8 & 24 & 10  & 7 & 24\\
				Eastern-Saxony-40-7 & 10 & 22 & 10 & 3 & 24\\
				Eastern-Saxony-40-8 & 8 & 22 & 20  & 3 & 24\\
				Eastern-Saxony-40-9 & 8 & 26 & 10  & 3 & 28\\
			\end{tabular}
	\end{table}
	\begin{table}[H]
		\centering
		\caption{Computational results for the robust test case \textit{Eastern-Saxony-40-2h}}
		\label{tab:ComputationalResultsRobust40}
			\begin{tabular}{lllllll}
				\textbf{Coverage}    & \textbf{Active}        & \textbf{Infrastructure}    & \textbf{Arcs}  & \textbf{Links} & \textbf{Runtime [s]}   & \textbf{Gap [\%]}  \\
				\textbf{Share [\%]}  & \textbf{Scenarios [\%]}& \textbf{Costs}             &       &   &               &           \\
				\hline
				10 & 10 & 847200 & 50 & 39 & 7200 & 12.78 \\
				20 & 20 & 881400 & 52 & 43 & 7200 & 15.30 \\
				30 & 30 & 1106100 & 64 & 49 & 7200 & 29.62 \\
				40 & 40 & 1025700 & 62 & 54 & 7200 & 25.26 \\
				50 & 50 & 888900 & 56 & 46 & 7200 & 8.17 \\
				60 & 70 & 922400 & 57 & 50 & 7200 & 10.00 \\
				70 & 70 & 972300 & 58 & 48 & 7200 & 14.09 \\
				80 & 80 & 1068800 & 62 & 51 & 7200 & 14.49 \\
				90 & 90 & 1176200 & 64 & 50 & 7200 & 14.99 \\
				100 & 100 & 1215800 & 67 & 49 & 7200 & 4.80 \\
			\end{tabular}
	\end{table}
	
	\section{Summary and Outlook}\label{sec:SummaryOutlook}
	
	This paper presents a formulation for a railway network design problem based on strategic timetables. By extracting the main properties of the timetable, which are the set of trains and their timing relationships, and using them as input for the network design problem, it is possible to tailor the network directly to the needs of the timetable. The infrastructure model allows for single- and multi-track railway lines and also considers which links in a node need to be included in the network solution. The capacity consumption on the railway lines is estimated through the use of train-type and train-sequence-dependent minimal headway times. Network capacity can be provided by including additional tracks or reducing travel or headway times. The initial strategic timetable is first reduced to an operational concept by using only the key properties and later refined during the optimization by assigning departure and arrival times to every train and every arc it is using.
	
	The model has been expanded to incorporate uncertainty within the strategic timetables. Uncertain timetables can be represented by a timetable family containing several discrete timetable scenarios and the robust model allows the calculation of networks which satisfy a given share of the input scenarios. Within one scenario, uncertain freight train demand can be modelled through optional trains. Using penalty coefficients, we can model different use cases for optional trains: very low penalties will make sure that optional trains will only be included in the solution if they can be added without additional infrastructure, while higher penalties induce a trade-off between penalties and additional infrastructure costs.
	
	Both the deterministic and the robust optimization models have been implemented using Python and solved by Gurobi. They have been tested on a case study using data derived from the current German strategic timetable concept, the \textit{Deutschlandtakt}. The computational results have shown that the optimization models work as intended, however, it has become apparent that the computational times are still too large to be able to optimize realistically sized instances. This is especially true for the robust multi-scenario case. 
	
	Because of that, improvements in computational times will be one major focus of future research. Both the deterministic and the robust formulation will be addressed, and decomposition algorithms as well as heuristic algorithms will be developed and evaluated. Besides, the quality of the resulting networks will be further evaluated by calculating the remaining capacity and conducting a sensitivity analysis to estimate the influence of input parameters on the resulting network and the objective value. The sensitivity analysis will focus on the cost ratios between penalties and infrastructure costs and between building and reduction costs.
	
	\section*{Acknowledgements}{This research has been funded by the Deutsche Forschungsgemeinschaft (DFG, German Research Foundation) - Project ID 432300662 and GRK 2236}
	
\printbibliography
	
\end{document}